\documentclass[preprint,fleqn,11pt]{elsarticle}
\usepackage{amssymb}
\biboptions{comma,square}

\journal{Applied Mathematics and Coputation}
\usepackage{amssymb,amsmath}
\usepackage{graphics}
\usepackage{graphicx}

\newcommand{\dl}[1]{{\bf Theorem{#1.}}}
\newcommand{\yl}[1]{{\bf Lemma{#1.}}}
\newcommand{\re}[1]{{\bf Remark{#1.}}}

\newcommand{\de}[1]{{\bf Definition{#1}}}
\newcommand{\zb}{\textbf{\qquad$\Box$}}
\newcommand{\la}[1]{\label{#1}}
\newcommand{\rf}[1]{(\ref{#1})}

\newcommand{\zm}{{\bf Proof.}}

\begin{document}
\thispagestyle{empty} \setcounter{page}{1}

\begin{frontmatter}

\title{A weak finite element method for elliptic problems in one space dimension
}

\author{Tie Zhang\footnote{Corresponding author at: Department of
Mathematics, Northeastern University, Shenyang 110004, China. {\em
E-mail address} : ztmath@163.com(T. Zhang). Tel \& Fax:
+86-024-83680949. },\quad Lixin Tang}

\address{Department of Mathematics and the State Key Laboratory of Synthetical Automation for Process
Industries, Northeastern University, Shenyang 110004, China}

\begin{abstract}
We present a weak finite element method for elliptic problems in
one space dimension. Our analysis shows that this method has more
advantages than the known weak Galerkin method proposed for
multi-dimensional problems, for example, it has higher accuracy
and the derived discrete equations can be solved locally, element
by element. We derive the optimal error estimates in the discrete
$H^1$-norm, the $L_2$-norm and the $L_\infty$-norm, respectively.
Moreover, some superconvergence results are also given. Finally,
numerical examples are provided to illustrate our theoretical
analysis.
\end{abstract}

\end{frontmatter}
\section{Introduction}
\setcounter{section}{1}\setcounter{equation}{0} Recently, the weak
Galerkin finite element method attracts much attention in the
field of numerical partial differential equations
\cite{Wang,Wang2,Lin,Lin1,Lin2,Lin3,Chen,Li,Har}. This method is
presented originally by Wang and Ye for solving elliptic problem
in multi-dimensional domain \cite{Wang}. Since then, some modified
weak Galerkin methods have also been studied, for example, see
\cite{Lin4,Gao,Wang0,Yang}. The weak Galerkin method can be
considered as an extension of the standard finite element method
where classical derivatives are replaced in the variational
equation by the weak derivatives defined on weak finite element 
functions. The main feature of this method is that it allows the
use of totally discontinuous finite element function and the trace
of finite element function on element boundary may be independent with
its value in the interior of element. This feature makes this
method possess the advantage of the usual discontinuous Galerkin
(DG) finite element method \cite{Arn,Cock,Zhang} and it has higher
flexibility than the DG method. The readers are referred to
articles \cite{Wang2,Lin,Cock} for more detailed explanation of
this method and its relation with other finite element methods.

In this paper, we present a weak finite element method for general
second order elliptic problem in one space dimension:
\begin{eqnarray}
\left\{
\begin{array}{lll}
-(a_2(x)u')'+a_1(x)u'+a_0(x)u=f(x),\,x\in (a,b),\\
u(a)=0, \; u'(b)=0,
\end{array}\label{1.1}
 \right.
\end{eqnarray}
where $a_2(x)\geq a_{min}>0,\,a_0(x)\geq 0$.

We first define the weak derivative and discrete weak
derivative on discontinuous function in one dimensional domain.
Then, we construct the weak finite element space $S_h$ and use it
to give the weak finite element approximation to problem \rf{1.1}.
Though, in some aspects, our method is similar to the original
weak Galerkin finite element method proposed for multi-dimensional
problem \cite{Wang}, it still has itself features. For example, we
impose the single value condition on space $S_h$ (see \rf{2.8} and
Remark 2.1), this condition can reduce the size of the finite
element discrete equations; Next, our space $S_h$ admits a weak
embedding inequality (see Lemma 3.2), which can be used to derive
the $L_\infty$-error estimate on mesh point set; Furthermore, the
discrete finite element system of equations derived from our
method can be solved locally, element by element, and this local
solvability is not feasible for the weak Galerkin method in
multi-dimensional space case. Except the usual optimal error
estimates in various norms, we also give some surperconvergence
results for the weak finite element solution. Numerical results
show that our method possesses very high computation accuracy. For
example, for finite element polynomial of order $k$, our
computation shows that the numerical convergence rates are at
least of order $k+2$ in the discrete $H^1$-norm, the $L_2$-norm
and the discrete $L_\infty$-norm. Our method also can be applied
to solve other partial differential equations in one space
dimension.

This paper is organized as follows. In Section 2, we introduce the
weak finite element method for the elliptic problem. In Section 3,
the stability of the weak finite element method is analyzed.
Section 4 is devoted to the optimal error estimate and
superconvergence estimate in various norms. In Section 5, the
local solvability of the weak finite element system of equations
is discussed and numerical experiments are provided to illustrate
our theoretical analysis.

Throughout this paper, we adopt the notations $H^{m}(I)$ to
indicate the usual Sobolev spaces on interval $I$ equipped with
the norm $\|\cdot\|_m=\|\cdot\|_{H^{m}(I)}$. The notations
$(\cdot,\cdot)$ and $\|\cdot\|$ denote the inner product and norm,
respectively, in the space $L_2(I)$. We will use letter $C$ to
represent a generic positive constant, independent of the mesh
size $h$.

\section{Problem and its weak finite element approximation}
\setcounter{section}{2}\setcounter{equation}{0} Consider elliptic
problem \rf{1.1}. Multiplying equation \rf{1.1} by the
transformation function
$$
\rho(x)=exp\,\Big(-\int_0^x\frac{a_1(x)}{a_2(x)}dx\Big),
$$
we see that problem \rf{1.1} can be transformed into the following
form:
$$
-(\rho\,a_2u')'+\rho\,a_0u=\rho\,f(x),\;x\in (a,b),\;
u(a)=0,\;u'(b)=0.
$$
Therefore, in what follows, we only consider elliptic problems in
the form:
\begin{eqnarray}
\left\{
\begin{array}{lll}
-(a_2(x)u')'+a_0(x)u=f(x),\,x\in (a,b),\\
u(a)=0, \; u'(b)=0,
\end{array}\label{2.1}
 \right.
\end{eqnarray}
where $a_2(x)\geq a_{min}>0,\,a_0(x)\geq 0$ and
$u'=\frac{du}{dx}$. We assume that $a_2(x)\in H^1(a,b),\,a_0(x)\in
L_\infty(a,b)$.

First, let us introduce the weak derivative concept. Let closed
interval $\bar{I}_a=[x_a, x_b]$ and its interior $I_a=(x_a,x_b)$.
A weak function on $\bar{I}_a$ refers to a function $v=\{
v^0,v^a,v^b\}$, $v^0=v|_{I_a}\in L_2(I_a)$, values $v^a=v(x_a)$
and $v^b=v(x_b)$ exist. Note that $v^a$ and $v^b$ may not be
necessarily the trace of $v^0$ at the interval endpoints $x_a$ and
$x_b$. Denote the weak function space by
$$
W(I_a)=\{ v=\{v^0,v^a,v^b\}: v^0\in L_2(I_a), |v^a|+|v^b|<\infty
\}.
$$
\de{ 2.1}\quad Let $v\in W(I_a)$. The weak derivative $d_wv$ of
$v$ is defined as a linear functional in the dual space
$H^{-1}(I_a)$ whose action on each $q\in H^1(I_a)$ is given by
\begin{equation}
<d_wv,q>\doteq-\int_{I_a}v^0q'dx+v^bq^b-v^aq^a,\;\forall\, q\in
H^1(I_a),\la{2.2}
\end{equation}
where $q^a=q(x_a), q^b=q(x_b)$.

Obviously, as a bounded linear functional on $H^1(I_a)$, $d_wv$ is
well defined for any $v\in W(I_a)$. Moreover, for $v\in H^1(I_a)$,
if we consider $v$ as a weak function with components
$v^0=v|_{I_a}$, $v^a=v(x_a)$ and $v^b=v(x_b)$, then by integration
by parts, we have for $q\in H^1(I_a)$ that
\begin{equation}
\int_{I_a}v'qdx=-\int_{I_a}vq'dx+v^bq^b-v^aq^a=-\int_{I_a}v^0q'dx+v^bq^b-v^aq^a,\la{2.3}
\end{equation}
which implies that $d_wv=v'$ is the usual derivative of function
$v$.

Next, we introduce the discrete weak derivative which is actually
used in our analysis. For nonnegative integer $r\geq 0$, let
$P_r(I_a)$ be the space composed of all polynomials on $I_a$ with
degree no more than $r$. Then, $P_r(I_a)$ is a subspace space of
$H^1(I_a)$.\\
\de{ 2.2}\quad For $v\in W(I_a)$, the discrete weak derivative
$d_{w,r}v\in P_r(I_a)$ is defined as the unique solution of the
following equation
\begin{equation}
\int_{I_a}d_{w,r}vqdx=-\int_{I_a}v^0q'dx+v^bq^b-v^aq^a,\;\forall\,
q\in P_r(I_a).\la{2.4}
\end{equation}

From \rf{2.2} and \rf{2.4}, we have
$$
<d_wv,q>=\int_{I_a}d_{w,r}vqdx,\;\forall\,q\in P_r(I_a).
$$
This shows that $d_{w,r}v$ is a discrete approximation of $d_wv$
in $P_r(I_a)$. In particular, if $v\in H^1(I_a)$, we have from
\rf{2.3} and \rf{2.4} that
$$
\int_{I_a}(d_{w,r}v-v')qdx=0,\;\forall\,q\in P_r(I_a).
$$
That is, $d_{w,r}v$ is the $L_2$ projection of $v'$ in $P_r(I_a)$
if $v\in H^1(I_a)$.

Now, we consider the weak finite element approximation of problem
\rf{2.1}. For interval $I=(a,b)$, let
$I_h:\;a=x_1<x_2<\cdots<x_{N-1}<x_N=b$ be a partition of $I$ with
elements $I_i=(x_i,x_{i+1}), i=1,\dots,N-1$. Denote the mesh size
by $h=\max h_i$, $h_i=x_{i+1}-x_i,\,i=1,\dots,N-1$. In the weak
finite element analysis, we will use the discrete weak function
space defined on partition $I_h$. Such space is denoted by
\begin{eqnarray}
W(I_h,k)&=&\{ v:\,v|_{I_i}\in W(I_i,k),\,i=1,2,\cdots, N-1\},\la{2.5}\\
W(I_i,k)&=&\{ v=\{v^0,v^i,v^{i+1}\}:\,v^0\in P_k(I_i),\,
|v^i|+|v^{i+1}|<\infty \}.\la{2.6}
\end{eqnarray}
Note that for a weak function $v\in W(I_i,k)$, the endpoint values
$v^i=v(x_i)$ and $v^{i+1}=v(x_{i+1})$ may be independent with the
interior value $v^0$. Recall the discrete weak derivative
definition \rf{2.4}, for $v\in W(I_i,k)$, its discrete weak
derivative $d_{w,r}v\in P_r(I_i)$ is given by the following
formula
\begin{equation}
\int_{I_i}d_{w,r}vqdx=-\int_{I_i}v^0q'dx+v^{i+1}q^{i+1}-v^iq^i,\;\forall\,q\in
P_r(I_i),\la{2.7}
\end{equation}
where $q^i=q(x_i), q^{i+1}=q(x_{i+1})$.

In our discussion, except for weak function $v=\{
v^0,v^i,v^{i+1}\}\in W(I_i,k)$, the endpoint values of a smooth
function $w$ on $I_i$ should be determined by its trace from the
interior of $I_i$. For example, for $w\in H^1(I_i)$,
$w^i=w(x_i)=\displaystyle{\lim_{x\rightarrow x_i}}w(x),\,x\in
I_i$.

Let $I_L=(x_{i-1},x_i)$ and $I_R=(x_i,x_{i+1})$ be two adjacent
elements with the common endpoint $x_i$, weak function
$v|_{\bar{I}_L}=\{ v^0_L,v^{i-1}_L, v^i_L\},\,v|_{\bar{I}_R}=\{
v^0_R,v^{i}_R, v^{i+1}_R\}$. We define the jump of weak function
$v$ at point $x_i$ by
$$
[v]_{x_i}=v^i_R-v^i_L,\;v\in W(I_h,k).
$$
Then, weak function $v$ is single value at point $x_i$ if and only
if $[v]_{x_i}=0$. Introduce the weak finite element space
\begin{equation}
S_h=\{ v:\,v\in W(I_h,k),\;
v^1=0,\,[v]_{x_i}=0,\,i=2,\dots,N-1\,\}.\la{2.8}
\end{equation}
Denote the discrete $L_2$ inner product and norm by
$$
(u,v)_h=\sum_{i=1}^{N-1}(u,v)_{I_i}=\sum_{i=1}^{N-1}\int_{I_i}u\,vdx,\;\;\;\;\|u\|_h^2=(u,u)_h.
$$
We now define the weak finite element approximation of problem
\rf{2.1} by finding $u_h\in S_h$ such that
\begin{equation}
(a_2d_{w,r}u_h,d_{w,r}v)_h+(a_0u^0_h,v^0)=(f,v^0),\;\forall\,v\in
S_h.\la{2.9}
\end{equation}
\re{ 2.1} The single value condition ($[v]_{x_i}=0$) has been
imposed on space $S_h$, it was not required in the original weak
Galerkin method \cite{Wang}. This condition can reduce the size of
discrete system of equations \rf{2.9}.

\section{The stability of weak finite element method}
\setcounter{section}{3}\setcounter{equation}{0} In this section,
we will show the stability of the weak finite element method and
give some lemmas which are very useful in our analysis.\\
\yl{ 3.1}\quad {\em Let $v\in W(I_i,k)$ and $r>k$. Then,
$d_{w,r}v=0$ if and only if $v=\{ v^0,v^{i},v^{i+1}\}$ is
constant on $\bar{I}_i$, that is, $v^0=v^{i}=v^{i+1}$ holds.}\\
\zm\quad First, let $v^0=v^{i}=v^{i+1}$. From \rf{2.7} we have
$$
\int_{I_i}d_{w,r}vqdx=v^0\Big(-\int_{I_i}q'dx+q^{i+1}-q^i\Big)=0,\,\forall\,q\in
P_r(I_i).
$$
This implies $d_{w,r}v=0$. Next, let  $d_{w,r}v=0$. Then we have
from \rf{2.7} that
\begin{equation}
-\int_{I_i}v^0q'dx+v^{i+1}q^{i+1}-v^iq^i=0,\,\forall\,q\in
P_r(I_i).\la{3.1}
\end{equation}
Let $\overline{v}=\frac{1}{h_i}\int_{I_i}v^0dx$ is the mean value
of $v$ on interval $I_i$. Consider the initial value problem:
\begin{eqnarray}
\left\{
\begin{array}{lll}
q'(x)=\overline{v}-v^0,\;x_{i}<x<x_{i+1},\;\;\\
q(x_i)=v^{i+1}-v^i.
\end{array}\label{3.2}
 \right.
\end{eqnarray}
Obviously, problem \rf{3.2} has a unique solution $q_1\in
P_r(I_i)$. By integrating \rf{3.2}, we obtain
$$
q_1^{i+1}-q^i_1=\int_{I_i}(\overline{v}-v^0)dx=0,\;\;\;q^{i+1}_1=q^i_1=v^{i+1}-v^i.
$$
Hence, taking $q=q_1$ in \rf{3.1}, we arrive at
$$
-\int_{I_i}v^0(\overline{v}-v^0)dx+(v^{i+1}-v^i)^2=\int_{I_i}(\overline{v}-v^0)^2dx
+(v^{i+1}-v^i)^2=0.
$$
This implies $v^0=\overline{v}$ and $v^i=v^{i+1}$. Substituting
this two equalities into \rf{3.1}, it yields
$$
(\overline{v}-v^i)(q^{i+1}-q^i)=0,\,\forall\,q\in P_r(I_i).
$$
Hence $\overline{v}=v^i$, so that $v^0=v^i=v^{i+1}$ holds. \zb

Lemma 3.1 shows that the discrete weak derivative $d_{w,r}v$
possesses the prominent feature of the classical derivative $v'$.
The following result is an analogy of the Sobolev embedding
theory in space $H^1_E(I)=\{\,v:\,v\in H^1(I),\,v(a)=0\}$. \\
\yl{ 3.2}\quad {\em Let $v\in S_h$ and $r>k$. Then, the following
weak embedding inequalities hold}
\begin{eqnarray}
&&|v^i|\leq
|x_i-a|^{\frac{1}{2}}\|d_{w,r}v\|_h,\;i=1,\dots,N,\,v\in
S_h,\la{3.3}\\
&&\|v^0\|\leq((b-a)+h)\|d_{w,r}v\|_h,\;v\in S_h.\la{3.4}
\end{eqnarray}
\zm\quad In definition \rf{2.7} of $d_{w,r}v$, taking $q=1$ we
have
\begin{equation}
\int_{I_i}d_{w,r}vdx=v^{i+1}-v^i\,.\la{3.5}
\end{equation}
Summing and using $v^1=0$ to obtain
\begin{eqnarray*}
v^{i+1}=\sum_{j=1}^{i}\int_{I_j}d_{w,r}vdx\leq\sum_{j=1}^i\sqrt{h_j}\|d_{w,r}v\|_{L_2(I_j)}
\leq |x_{i+1}-a|^{\frac{1}{2}}\|d_{w,r}v\|_h.
\end{eqnarray*}
This gives estimate \rf{3.3}. To prove \rf{3.4}, let $q_1\in
P_r(I_i)$ satisfies the initial problem:
\begin{eqnarray}
\left\{
\begin{array}{lll}
q'_1(x)=-v^0,\;x_{i}<x<x_{i+1},\;\;\\
q_1(x_i)=v^{i+1}-v^i.
\end{array}\label{3.6}
 \right.
\end{eqnarray}
Taking $q=q_1$ in \rf{2.7}, we obtain
\begin{equation}
\int_{I_i}|v^0|^2dx=\int_{I_i}d_{w,r}vq_1dx+v^iq^i_1-v^{i+1}q_1^{i+1}.\la{3.7}
\end{equation}
Integrating \rf{3.6}, it yields
\begin{equation}
q_1(x)=v^{i+1}-v^i-\int_{x_i}^xv^0dx,\;\;\;q_1^{i+1}=v^{i+1}-v^i-\int_{I_i}v^0dx.\la{3.8}
\end{equation}
Substituting \rf{3.8} into \rf{3.7} and using \rf{3.5}, we obtain
\begin{eqnarray*}
\|v^0\|^2_{L_2(I_i)}&=&\int_{I_i}d_{w,r}v
dx(v^{i+1}-v^i)-\int_{I_i}d_{w,r}v\int_{x_i}^xv^0(y)dydx\\
&&+v^i(v^{i+1}-v^i)-v^{i+1}\Big(v^{i+1}-v^i-\int_{I_i}v^0dx\Big)\\
&=&-\int_{I_i}d_{w,r}v\int_{x_i}^xv^0(y)dy
dx+v^{i+1}\int_{I_i}v^0dx.
\end{eqnarray*}
Hence, it follows from estimate \rf{3.3} and the Cauchy inequality
that
\begin{eqnarray*}
\|v^0\|^2\leq h\|d_{w,r}v\|_h\|v^0\|+(b-a)\|d_{w,r}v\|_h\|v^0\|.
\end{eqnarray*}
The proof is completed. \zb

Now we can prove the stability of weak finite element equation
\rf{2.9}.\\
\dl{ 3.3}\quad{\em Let $r>k$. Then problem \rf{2.9} has a unique
solution $u_h\in S_h$ and $u_h$ satisfies the stability estimate}
\begin{equation}
\|u_h^0\|+\|d_{w,r}u_h\|_h\leq
\frac{2((b-a)+1)^2}{a_{min}}\|f\|.\la{3.9}
\end{equation}
\zm\quad First, consider the stability. Taking $v_h=u_h$ in
\rf{2.9}, we have
$$
a_{min}\|d_{w,r}u_h\|_h^2\leq \|f\|\,\|u_h^0\|\,.
$$
Together with \rf{3.4}, estimate \rf{3.9} is derived. Next,
consider the unique existence. Since problem \rf{2.9} is a linear
system composed of $(k+2)\times (N-1)$ equations with $(k+2)\times
(N-1)$ unknowns, we only need to prove that $u_h=0$ if $f=0$. Let
$f=0$, then it follows from \rf{3.9} that
$\|d_{w,r}u_h\|_h=\|u_h^0\|=0$ holds. Therefore, from Lemma 3.1,
we can conclude that $u_h$ is piecewise constant on partition
$I_h$ so that $u_h^i=u_h^0=0$. \zb

\section{Error analysis}
\setcounter{section}{4}\setcounter{equation}{0} In this section,
we do the error analysis for the weak finite element method
\rf{2.9}. We will see that the weak finite element method
possesses the same or better theoretical convergence rate as that
of the conventional finite element method.

We first show the approximation property of the weak finite
element space $S_h$. In order to balance the approximation
accuracy between space $S_h$ and space $P_r(I_i)$ used for
$d_{w,r}v$, from mow on, we always set the index $r=k+1$ in the
definition of discrete weak derivative $d_{w,r}v$, see \rf{2.7}.

For $l\geq 0$, let $P_h^l$ is the local $L_2$ projection operator,
restricted on each element $I_i$, $P_h^l:\,u\in
L_2(I_i)\rightarrow P_h^lu\in P_l(I_i)$ such that
\begin{equation}
(u-P_h^lu,q)_{I_i}=0,\;\forall\,q\in
P_l(I_i),\,i=1,2,\dots,N-1.\la{4.1}
\end{equation}
By the Bramble-Hilbert lemma, it is easy to prove that (see
\cite{Zhang})
\begin{equation}
\|u-P_h^lu\|_{L_2(I_i)}\leq Ch_i^{s}\|u\|_{H^{s}(I_i)},\;0\leq
s\leq l+1.\la{4.2}
\end{equation}
We now define a projection operator $Q_h:\,u\in H^1(I)\rightarrow
Q_hu\in W(I,k)$ such that
\begin{equation}
Q_hu|_{\bar{I}_i}=\{Q^0_hu,(Q_hu)^i
,(Q_hu)^{i+1}\}\doteq\{P_h^ku,u(x_i),u(x_{i+1})\},\;i=1,\dots,N-1.\la{4.3}
\end{equation}
Obviously, $Q_hu\in S_h$ if $u\in H^1_E(I)$. From \rf{4.2}, we
have
\begin{equation}
\|Q_h^0u-u\|_{L_2(I_i)}=\|P_h^ku-u\|_{L_2(I_i)}\leq
Ch_i^s\|u\|_{H^s(I_i)},\;0\leq s\leq k+1.\la{4.4}
\end{equation}
Furthermore, since
\begin{eqnarray*}
\int_{I_i}d_{w,r}Q_huqdx&=&-\int_{I_i}Q_h^0uq'dx+(Q_hu)^{i+1}q^{i+1}-(Q_hu)^iq^i\\
&=&-\int_{I_i}u
q'dx+u^{i+1}q^{i+1}-u^iq^i=\int_{I_i}u'qdx,\;\;\forall\,q\in
P_r(I_i),
\end{eqnarray*}
hence $d_{w,r}Q_hu=P_h^ru'$ holds and (noting that $r=k+1$)
\begin{equation}
\|d_{w,r}Q_hu-u'\|_{L_2(I_i)}=\|P_h^ru'-u'\|_{L_2(I_i)}\leq
Ch_i^s\|u\|_{H^{s+1}(I_i)},\;0\leq s\leq k+2.\la{4.5}
\end{equation}
Estimates \rf{4.4} and \rf{4.5} show that $Q_hu\in S_h$ is a very
good approximation for function $u\in H^1_E(I)\bigcap H^m(I),
m\geq 1$.

In order to do the error analysis, we still need to construct
another special projection function.\\
\yl{ 4.1}\quad{\em For $u\in H^1(I)$, there exists a projection
function $\pi_hu\in H^1(I)$, restricted on element $I_i$,
$\pi_hu\in P_{k+1}(I_i)$ satisfies}
\begin{eqnarray}
&&((\pi_hu)',q)_{I_i}=(u',q)_{I_i},\;\forall\,q\in
P_k(I_i),\,i=1,\dots,N-1,\la{4.6}\\
&&\pi_hu(x_i)=u(x_i),\;i=1,\dots,N,\la{4.7}\\
&&\|u-\pi_hu\|_{L_2(I_i)}+h_i\|u-\pi_hu\|_{H^1(I_i)}\leq
Ch_i^{s+1}\|u\|_{H^{s+1}(I_i)},\,0\leq s\leq k+1.\la{4.8}
\end{eqnarray}
\zm\quad Let $u\in H^1(I)$. For any given element $I_i$, let
$\pi_h^{(i)}u\in P_{k+1}(I_i)$ be the unique solution of the
initial problem:
\begin{eqnarray}
\left\{
\begin{array}{lll}
(\pi_h^{(i)}u)'(x)=P_h^ku',\;x_{i}<x<x_{i+1},\;\;\\
\pi_h^{(i)}u(x_i)=u(x_i).
\end{array}\label{4.9}
 \right.
\end{eqnarray}
Then, by the property of operator $P_h^k$, we obtain
\begin{eqnarray}
&&((\pi_h^{(i)}u)',q)_{I_i}=(u',q)_{I_i},\;\forall\,q\in
P_k(I_i).\la{4.10}\\
&&\|u'-(\pi_h^{(i)}u)'\|_{L_2(I_i)}\leq
Ch_i^{s}\|u\|_{H^{s+1}(I_i)},\,0\leq s\leq k+1.\la{4.11}
\end{eqnarray}
Since
$$
(\pi_h^{(i)}u-u)(x)=\int_{x_i}^x(\pi_h^{(i)}u-u)'(x)dx,\;x\in I_i,
$$
hence, it follows from \rf{4.10} and the Cauchy inequality that
\begin{equation}
\pi_h^{(i)}u(x_{i+1})=u(x_{i+1}),\;\;\;\|u-\pi_h^{(i)}u\|_{L_2(I_i)}\leq
h_i\|u'-(\pi_h^{(i)}u)'\|_{L_2(I_i)}.\la{4.12}
\end{equation}
Now, we set $\pi_hu|_{I_i}=\pi_h^{(i)}u$ for $1\leq i\leq N-1$,
then conclusions \rf{4.6}$\sim$\rf{4.8} can be derived by using
\rf{4.10}$\sim$\rf{4.12}. Furthermore, since
$\pi_h^{(i)}u(x_{i+1})=u(x_{i+1})=\pi_h^{(i+1)}u(x_{i+1})$, this
shows that $\pi_hu$ is continuous across junction point $x_{i+1}$,
so $\pi_hu\in H^1(I)$ holds. \zb\\
\yl{ 4.2}\quad {\em Let $u\in H^1_E(I)\bigcap H^2(I)$ be the
solution of problem \rf{2.1}. Then, $u$ satisfies the following
equation}
\begin{equation}
(\pi_h(a_2u'),d_{w,r}v)_h+(a_0u,v^0)=(f,v^0),\;\forall\,v\in
S_h.\la{4.13}
\end{equation}
\zm\quad By \rf{2.7} and Lemma 4.1, we have for $v\in S_h$ that
\begin{eqnarray*}
(\pi_h(a_2u'),d_{w,r}v)_{I_i}&=&-((\pi_h(a_2u'))',v^0)_{I_i}+(\pi_h(a_2u'))^{i+1}v^{i+1}-(\pi_h(a_2u'))^iv^{i}\\
&=&-((a_2u')',v^0)_{I_i}+(a_2u')^{i+1}v^{i+1}-(a_2u')^iv^{i}.
\end{eqnarray*}
Summing and noting that $v^1=0$ and $u'(x_N)=0$, it yields
$$
(\pi_h(a_2u'),d_{w,r}v)_h=-((a_2u')',v^0)=-(a_0u,v^0)+(f,v^0).
$$
Hence, equation \rf{4.13} holds. \zb\\
\dl{ 4.3}\quad{\em Let $u$ and $u_h$ be the solutions of problems
\rf{2.1} and \rf{2.9}, respectively, $u\in H^1_E(I)\bigcap H^2(I)$
and $r=k+1$. Then we have}
\begin{equation}
a_{min}\|d_{w,r}Q_hu-d_{w,r}u_h\|_h\leq
\|a_2d_{w,r}Q_hu-\pi_h(a_2u')\|_h+((b-a)+1)\|a_0(Q^0_hu-u)\|.\la{4.14}
\end{equation}
\zm\quad From Lemma 4.2, we have
\begin{eqnarray*}
&&(a_2d_{w,r}Q_hu,d_{w,r}v)_h+(a_0Q^0_hu,v^0)\\
&=&(f,v^0)+(a_2d_{w,r}Q_hu-\pi_h(a_2u'),d_{w,r}v)_h+(a_0(Q^0_hu-u),v^0).
\end{eqnarray*}
Combining this with equation \rf{2.9}, we obtain the error
equation
\begin{eqnarray}
&&(a_2d_{w,r}(Q_hu-u_h),d_{w,r}v)_h+(a_0(Q_h^0u-u_h^0),v^0)\nonumber\\
&=&(a_2d_{w,r}Q_hu-\pi_h(a_2u'),d_{w,r}v)_h+(a_0(Q^0_hu-u),v^0),\;v\in
S_h.\la{4.15a}
\end{eqnarray}
Taking $v=Q_hu-u_h\in S_h$ and using the weak embedding equality
\rf{3.4}, we arrive at the conclusion of Theorem 4.3. \zb

By means of Theorem 4.3, we can derive the following error
estimates.\\
\dl{ 4.4}\quad{\em Let $u$ and $u_h$ be the solutions of problems
\rf{2.1} and \rf{2.9}, respectively, $u\in H^1_E(I)\bigcap
H^{2+s}(I), a_2\in H^{1+s}(I), s\geq 0$, and $r=k+1$. Then we have
\begin{eqnarray}
\|d_{w,r}u_h-u'\|_h\leq Ch^{s+1}\|u\|_{s+2},\;0\leq s\leq k,\la{4.15}\\
\max_{1\leq i\leq N}|u_h^i-u(x_i)|\leq Ch^{s+1}\|u\|_{s+2},\;0\leq
s\leq k.\la{4.16}
\end{eqnarray}
Furthermore, if $a_0(x)=0$ and $u$ is smooth enough, then we have
the superconvergence estimates}
\begin{eqnarray}
\|d_{w,r}u_h-u'\|_h\leq Ch^{k+2}\|u\|_{k+3},\;k\geq 0,\la{4.17}\\
\max_{1\leq i\leq N}|u_h^i-u(x_i)|\leq Ch^{k+2}\|u\|_{k+3},\;k\geq
0.\la{4.18}
\end{eqnarray}
\zm\quad By the triangle inequality, we have
\begin{eqnarray*}
&&\|d_{w,r}u_h-u'\|_h\leq
\|d_{w,r}u_h-d_{w,r}Q_hu\|_h+\|d_{w,r}Q_hu-u'\|_h,\\
&&\|a_2d_{w,r}Q_hu-\pi_h(a_2u')\|_h\leq
\|a_2(d_{w,r}Q_hu-u')\|_h+\|a_2u'-\pi_h(a_2u')\|.
\end{eqnarray*}
Together with Theorem 4.3, it yields
\begin{equation}
\|d_{w,r}u_h-u'\|_h\leq
C\big(\,\|d_{w,r}Q_hu-u'\|_h+\|a_2u'-\pi_h(a_2u')\|+\|a_0(Q^0_hu-u)\|\,\big).\la{4.19a}
\end{equation}
Then, estimate \rf{4.15} follows from the approximation properties
\rf{4.4}, \rf{4.5} and \rf{4.8}. Furthermore, by the weak
embedding inequality \rf{3.3}, we have
$$
|u(x_i)-u_h^i|=|(Q_hu)^i-u_h^i|\leq
(x_i-a)^{\frac{1}{2}}\|d_{w,r}Q_hu-d_{w,r}u_h\|_h.
$$
Hence, we can obtain estimate \rf{4.16} by using Theorem 4.3 and
the approximations properties. The superconvergence estimates
\rf{4.17}-\rf{4.18} can be derived by a similar argument, noting
that $\|a_0(Q^0_hu-u)\|=0$ in \rf{4.14} and \rf{4.19a} if $a_0=0$.
\zb

From Theorem 4.3 and the weak embedding inequality, we immediately
obtain
\begin{equation}
\|Q_h^0u-u_h^0\|\leq C\|d_{w,r}Q_hu-d_{w,r}u_h\|_h\leq
Ch^{s+1}\|u\|_{s+2},\;0\leq s\leq k\,,\la{4.19}
\end{equation}
which results in the $L_2$ error estimate
$$
\|u-u_h^0\|\leq \|u-Q^0_hu\|+\|Q^0_hu-u_h^0\|\leq
Ch^{s+1}\|u\|_{s+2},\;0\leq s\leq k.
$$

Below we give a superclose estimate for error $Q^0_hu-u_h^0$. To
this end, we introduce the auxiliary problem: Find $w\in
H^1_E(I)\bigcap H^2(I)$ such that
\begin{eqnarray}
\left\{
\begin{array}{lll}
-(a_2(x)w')'+a_0(x)w=Q_h^0u-u^0_h,\,x\in (a,b),\\
w(a)=0, \; w'(b)=0,\;\;\|w\|_2\leq C\|Q^0_hu-u^0_h\|.
\end{array}\label{4.20}
 \right.
\end{eqnarray}
From Lemma 4.2, we know that $w$ satisfies equation:
\begin{equation}
(\pi_h(a_2w'),d_{w,r}v)_h+(a_0w,v^0)=(Q^0_hu-u^0_h,v^0),\;\forall\,v\in
S_h.\la{4.21}
\end{equation}
\dl{ 4.5}\quad{\em Let $u$ and $u_h$ be the solutions of problems
\rf{2.1} and \rf{2.9}, respectively, $u\in H^1_E(I)\bigcap
H^{2+s}(I), a_2\in H^{1+s}(I),\,a_0\in H^1(I), s\geq 0$, and
$r=k+1$. Then we have the following superclose estimate}
\begin{eqnarray}
\|Q_h^0u-u_h^0\|\leq Ch^{s+2}\|u\|_{s+2},\;0\leq s\leq k.\la{4.22}
\end{eqnarray}
\zm\quad Taking $v=Q_hu-u_h$ in \rf{4.21} and using error equation
\rf{4.15a}, we have
\begin{eqnarray}
&&\|Q_h^0u-u_h^0\|^2\nonumber\\
&=&(d_{w,r}(Q_hu-u_h),\pi_h(a_2w'))_h+(a_0(Q_h^0u-u_h^0),w)\nonumber\\
&=&(d_{w,r}(Q_hu-u_h),\pi_h(a_2w')-a_2d_{w,r}Q_hw)_h+(a_0(Q_h^0u-u_h^0),w-Q_h^0w)\nonumber\\
&&+(a_2d_{w,r}(Q_hu-u_h),d_{w,r}Q_hw)_h+(a_0(Q_h^0u-u_h^0),Q_h^0w)\nonumber\\
&=&(d_{w,r}(Q_hu-u_h),\pi_h(a_2w')-a_2d_{w,r}Q_hw)_h+(a_0(Q_h^0u-u_h^0),w-Q_h^0w)\nonumber\\
&&+(a_2d_{w,r}Q_hu-\pi_h(a_2u'),d_{w,r}Q_hw)_h+(a_0(Q_h^0u-u),Q_h^0w)\nonumber\\
&=&\big\{(d_{w,r}(Q_hu-u_h),\pi_h(a_2w')-a_2d_{w,r}Q_hw)_h+(a_0(Q_h^0u-u_h^0),w-Q_h^0w)\big\}\nonumber\\
&&+\big\{(a_2d_{w,r}Q_hu-\pi_h(a_2u'),d_{w,r}Q_hw-w')_h+(a_0(Q_h^0u-u),Q_h^0w-w)\big\}\nonumber\\
&&+\big\{(a_2d_{w,r}Q_hu-\pi_h(a_2u'),w')_h+(a_0(Q_h^0u-u),w)\big\}\nonumber\\
&=&E_1+E_2+E_3.\la{4.24}
\end{eqnarray}
Below we estimate $E_1\sim E_3$. Using \rf{4.19} and the
approximation properties of operators $Q_h$ and $\pi_h$, we have
\begin{eqnarray*}
E_1&=&(d_{w,r}(Q_hu-u_h),\pi_h(a_2w')-a_2d_{w,r}Q_hw)_h+(a_0(Q_h^0u-u_h^0),w-Q_h^0w)\\
&\leq&
C\|d_{w,r}(Q_hu-u_h)\|_h\big(\|\pi_h(a_2w')-a_2w'+a_2w'-a_2d_{w,r}Q_hw\|_h+\|w-Q_h^0w\|\big)\\
&\leq&Ch^{s+2}\|u\|_{s+2}\|w\|_2.\\
E_2&=&(a_2d_{w,r}Q_hu-\pi_h(a_2u'),d_{w,r}Q_hw-w')_h+(a_0(Q_h^0u-u),Q_h^0w-w)\\
&\leq&
Ch\big(\|a_2d_{w,r}Q_hu-a_2u'+a_2u'-\pi_h(a_2u')\|_h\big)\|w\|_2+Ch^{s+2}\|u\|_{s+1}\|w\|_1\\
&\leq&Ch^{s+2}\|u\|_{s+2}\|w\|_2.
\end{eqnarray*}
Next, we write
\begin{eqnarray*}
E_3&=&(a_2d_{w,r}Q_hu-\pi_h(a_2u'),w')_h+(a_0(Q_h^0u-u),w)\\
&=&(a_2d_{w,r}Q_hu-a_2u',w')_h+(a_2u'-\pi_h(a_2u'),w')_h+(a_0(Q_h^0u-u),w)\\
&=&E_{31}+E_{32}+E_{33}.
\end{eqnarray*}
Since $d_{w,r}Q_hu=P_h^ru',\,Q_h^0u=P_h^ku$, then we have
\begin{eqnarray*}
E_{31}+E_{33}&=&(P_h^ru'-u',a_2w'-P_h^k(a_2w'))_h+(P_h^ku-u,a_0w-P_h^k(a_0w))\\
&\leq&Ch^{s+2}\|u\|_{s+2}\|w\|_2.
\end{eqnarray*}
Furthermore, from Lemma 4.1 and integration by parts, we also
obtain
\begin{eqnarray*}
E_{32}&=&-\sum_{i=1}^{N-1}((a_2u'-\pi_h(a_2u'))',w)_{I_i}=-\sum_{i=1}^{N-1}((a_2u')'-P_h^k(a_2u')',w-P_h^kw)_{I_i}\\
&\leq&Ch^{s+2}\|u\|_{s+2}\|w\|_2.
\end{eqnarray*}
Hence, we have that $E_3\leq Ch^{s+2}\|u\|_{s+2}\|w\|_2$. The
proof is completed by substituting estimates $E_1\sim E_3$ into
\rf{4.24}, noting that $\|w\|_2\leq C\|Q_h^0u-u_h^0\|$.\zb

From Theorem 4.5 and the triangle inequality, we immediately
obtain the following optimal $L_2$-norm error estimate
\begin{equation}
\|u-u_h^0\|\leq Ch^{k+1}\|u\|_{k+1},\;k\geq 1.
\end{equation}

In order to derive the optimal $L_\infty$-error estimate, we need
to strengthen the partition condition. Partition $I_h$ is called
quasi-uniform if there exists a positive constant $\sigma$ such
that
$$
h/h_i\leq\sigma,\, i=1,\cdots,N.
$$
This condition assures that the inverse inequality holds in space
$S_h$.\\
\dl{ 4.6}\quad{\em Assume that partition $I_h$ is quasi-uniform,
and $u$ and $u_h$ are the solution of problems \rf{2.1} and
\rf{2.9}, respectively, and conditions in Theorem 4.5 hold. Then,
we have}
\begin{equation}
\|u-u_h^0\|_{L_\infty(I)}\leq Ch^{s+1}\|u\|_{s+2}, \,0\leq s\leq
k.
\end{equation}
\zm\quad From Theorem 4.5 and the finite element inverse
inequality, we have that
$$
\|Q_h^0u-u_h^0\|_{L_\infty(I)}\leq
Ch^{-\frac{1}{2}}\|Q_h^0u-u_h^0\|\leq
Ch^{s+\frac{3}{2}}\|u\|_{s+2}.
$$
Hence, by using the approximation property of $Q_h^0u=P_h^ku$, we
obtain
\begin{eqnarray*}
\|u-u_h^0\|_{L_\infty(I)}&\leq&
\|u-Q_h^0u\|_{L_\infty(I)}+\|Q_h^0u-u_h^0\|_{L_\infty(I)}\\
&\leq& Ch^{s+1}(\|u\|_{s+1,\infty}+\|u\|_{s+2})\leq
Ch^{s+1}\|u\|_{s+2},
\end{eqnarray*}
where we have used the Sobolev embedding inequality. \zb
\section{The local solvability and numerical example}
\setcounter{section}{5} \setcounter{equation}{0} In this section,
we discuss how to solve the discrete system of equations \rf{2.9}.
We will design a local solver so that this linear system can be
solved locally, element by element. Then, we provide some
numerical examples to illustrate our theoretical analysis.

\subsection{The local solvability of the weak finite element
equation}

Consider the weak finite element equation: (see \rf{2.9}):
\begin{equation}
(a_2d_{w,r}u_h,d_{w,r}v)_h+(a_0u^0_h,v^0)=(f,v^0),\;\forall\,v\in
S_h.\la{5.1}
\end{equation}
In order to form the discrete linear system of equations \rf{5.1},
we introduce the basis functions of space $W(I_i,k)$ or $S_h$. Let
weak basis functions
$\psi_{j}(x)=\{\psi^0_j,\psi^i_j,\psi^{i+1}_j\}=\{x^{j-1},0,0\},\,j=1,\dots,k+1$,
and further let $\delta_i(x)$ be the node basis function, that is,
$\delta_i(x_i)=1, \delta_i(x)=0, x\neq x_i$. Then, we have
$W(I_i,k)=span\{\psi_1(x),\dots,\psi_{k+1}(x),\delta_i(x),$
$\delta_{i+1}(x)\}$, and for any $v\in S_h$, restricted on $I_i$,
$v=\{v^0,v^i,v^{i+1}\}$ can be written as
\begin{eqnarray*}
v(x)=\sum_{j=1}^{k+1}c_j\psi_j(x)+v^i\delta_i(x)+v^{i+1}\delta_{i+1}(x),\;x_i\leq
x\leq x_{i+1}.
\end{eqnarray*}
For $v\in S_h$, by the definition \rf{2.7} of discrete weak
derivative, we see that the support set of $d_{w,r}\psi_j(x)$ is
in $I_i$ and the support set of $d_{w,r}\delta_i(x)$ is in
$\bigcup I_j$, where $\overline{I}_j\bigcap x_i\neq\varnothing$.
Then, equation \rf{5.1} is equivalent to the following system of
equations
\begin{eqnarray}
&&(a_2d_{w,r}u_h,d_{w,r}v)_{I_i}+(a_0u^0_h,v^0)_{I_i}=(f,v^0)_{I_i},\,v=\{\psi_j\},
\,i=1,\dots, N-1,  \la{5.2}\\
&&(a_2d_{w,r}u_h,d_{w,r}v)_{I_{i}\cup
I_{i+1}}=0,\,v=\delta_{i+1},\,i=1,\dots,N-2,\la{5.3}\\
&&(a_2d_{w,r}u_h,d_{w,r}v)_{I_{N-1}}=0,\,v=\delta_{N}.\la{5.4}
\end{eqnarray}
Equations \rf{5.2}$\sim$\rf{5.4} form a linear system composed of
$(k+2)(N-1)$ equations with $(k+2)(N-1)$ unknowns. To solve this
system, we need to design a solver for the discrete weak
derivative $d_{w,r}v$ or $d_{w,r}u_h$. According to \rf{2.7}, for
given $v\in W(I_i,k)$, $d_{w,r}v\in P_r(I_i)$ can be computed by
the following formula
\begin{equation}
M_id_{w,r}V=A_iV^0+B_iV^i, \la{5.5}
\end{equation}
where $d_{w,r}V$ and $V^0$ are the vectors associated with
functions $d_{w,r}v\in P_r(I_i)$ and $v^0\in P_k(I_i)$,
respectively, and $V^i=(v^i, v^{i+1})^T$. The matrixes in \rf{5.5}
are as follows
\begin{eqnarray*}
&&M_i=(m_{st})_{(r+1)\times (r+1)}, \,\,A_i=(a_{st})_{(r+1)\times
(k+1)},\,\,B_i=(b_{st})_{(r+1)\times 2},\\
&&m_{st}=(x^{s-1},x^{t-1})_{I_i},\,\,a_{st}=-(d_x(x^{s-1}),x^{t-1})_{I_i},\,
\,b_{s1}=-x_i^{s-1},\,b_{s2}=x_{i+1}^{s-1}.
\end{eqnarray*}
Now, linear system of equations \rf{5.2}$\sim$\rf{5.4} can be
solved in the following two ways.

{\bf Method One}. We first use formula \rf{5.5} to derive the
linear representation $d_{w,r}u_h(I_i)$
$=L(u_h^0(I_i),u^i_h,u^{i+1}_h)$. Then, by substituting
$d_{w,r}u_h(I_i)$ into equations \rf{5.2}$\sim$\rf{5.4}, we can
obtain a linear system of equations that only concerns unknowns
$\{u_h^0(I_i),u_h^i,u_h^{i+1}\}, i=1,\dots,N-1$. Now, this linear
system can be solved by using a proper linear solver, in which
$d_{w,r}v(I_i)$ is computed by formula \rf{5.5}.

{\bf Method Two}. We observe that the unknowns in equations
\rf{5.2}$\sim$\rf{5.4} are coupled only by equation \rf{5.3} which
concerning unknowns on two adjacent elements. If we can
independently solve the unknowns on some single element, then we
are able to uncouple this simultaneous equations and solve the
whole linear system of equations \rf{5.2}$\sim$\rf{5.4} locally,
element by element. To this end, integrating equation \rf{2.1}, we
find that the exact solution $u$ satisfies
$$
u(x_N)-u(x_{N-1})+\int_{I_{N-1}}\frac{1}{a_2(x)}\widetilde{u}(x)dx=
\int_{I_{N-1}}\frac{1}{a_2(x)}\widetilde{f}(x)dx,
$$
where
$$
\widetilde{u}(x)=\int^{x_N}_x
a_0(y)u(y)dy,\;\;\widetilde{f}(x)=\int^{x_N}_xf(y)dy.
$$
This provides an additional equation for $u_h$ on the last
element. Now, we can solve linear system of equations
\rf{5.2}$\sim$\rf{5.4} locally in the following procedure.

First, on element $I_{N-1}$, solve $u_h=(u_h^0,u_h^{N-1},u_h^N)$
by the equations:
\begin{eqnarray*}
&&(a_2d_{w,r}u_h,d_{w,r}v)_{I_{N-1}}+(a_0u^0_h,v^0)_{I_{N-1}}=(f,v^0)_{I_{N-1}},\,v=\psi_j,\,
j=1,\dots,k+1,\\
&&(a_2d_{w,r}u_h,d_{w,r}v)_{I_{N-1}}=0,\,v=\delta_{N},\\
&&u_h^N-u_h^{N-1}+\int_{I_{N-1}}\frac{1}{a_2(x)}\widetilde{u}_h^0(x)dx=
\int_{I_{N-1}}\frac{1}{a_2(x)}\widetilde{f}(x)dx.
\end{eqnarray*}
Then, on each element $I_i$, solve $u_h=(u_h^0,u_h^{i},u_h^{i+1})$
by the following equations in the order of $i=N-2,\dots,1$,
\begin{eqnarray*}
&&(a_2d_{w,r}u_h,d_{w,r}v)_{I_{i}}+(a_0u^0_h,v^0)_{I_i}=(f,v^0)_{I_i},\,v=\psi_j,
\,j=1,\dots,k+1,\\
&&(a_2d_{w,r}u_h,d_{w,r}v)_{I_{i}}=-(a_2d_{w,r}u_h,d_{w,r}v)_{I_{i+1}},\,v=\delta_{i+1},\;\;[u_h]_{x_{i+1}}=0.
\end{eqnarray*}
In the above computation procedure, $d_{w,r}v$ and $d_{w,r}u_h$
are still determined by formula \rf{5.5}. It is easy to see that
Method Two is more economical than Method One.

\subsection{Numerical example}

Let us consider problem \rf{2.1} with the following data:
\begin{equation}
u(x)=2(1-x)\sin(\pi x),\; a_2(x)=1+x^2,\;a_0(x)=\sin (\pi
x),\la{5.6}
\end{equation}
and $(a,b)=(0, 1)$, the corresponding source term
$f=-(a_2u')'+a_0u$.

In the numerical experiments, we always partition the interval
$I=(0,1)$ uniformly with the mesh size $h=1/N$. We examine the
computation error in the discrete $H^1$-norm, the $L_2$-norm and
the $L_\infty$-norm on the mesh point set. The numerical
convergence rate is computed by using the formula
$r=\ln(e_h/e_{\frac{h}{2}})/\ln 2$, where $e_h$ is the computation
error. Table 5.1$\sim$Table 5.3 give the numerical results with
finite element polynomials of order $k=0,1,2$, in sequence. We
observe that the errors vanish very quickly and the convergence
rates are at least one order higher than that theoretically
predicted, i.e., the superconvergence results are obtained even
$a_0(x)\neq 0$. When taking $a_0(x)=0$, we obtain the same
superconvergence rate as that in case of $a_0(x)\neq 0$. We
further examine problem \rf{2.1} with different test solutions and
data, the convergence rates still remain unchanged. In conclusion,
this weak finite element method is a high accuracy numerical
method in both theory and experiment.

\begin{center}
{\small{\bf Table 5.1}\quad History of convergence for $k=0$}\\[0.5\baselineskip]
\renewcommand\arraystretch{1.4}
\arrayrulewidth 0.5pt \small
\begin{tabular}{cccc}
\hline
&$\|d_{w,r}u_h-u'\|_h$&\quad$\|u-u_h^0\|$&\quad$\max|u_h^i-u(x_i)|$\\
mesh $h$&error \quad\quad rate&\quad error \qquad\quad rate&\quad error \qquad\quad rate\\
\hline
1/4&0.2281\quad\ \ \ \ \ \ \ - &\quad 0.0501 \qquad\ \ \ - &\quad 0.1221 \qquad \ \ \ \ \ \ \ \ \ -\\
1/8&0.0579\quad 1.9769 &\qquad0.0131\qquad 1.9361         &\quad 0.0302\qquad\ \ \ 2.0162\\
1/16&0.0145\quad 1.9942 &\qquad0.0033\qquad 1.9836        &\quad 0.0075\qquad\ \ \ 2.0039\\
1/32&0.0036\quad 1.9986 &\qquad 0.0008\qquad 1.9959       &\quad 0.0019\qquad\ \ \ 2.0010\\
1/64&0.0009\quad 1.9996 &\qquad 0.0002\qquad 1.9990       &\quad 0.0005\qquad\ \ \ 2.0002\\
1/128&0.0002\quad1.9999 &\qquad0.0001\qquad 1.9997        &\quad 0.0001\qquad\ \ \ 2.0001\\
 \hline
\end{tabular}
\end{center}
\vspace{0.2cm}

\begin{center}
{\small{\bf Table 5.2}\quad History of convergence for $k=1$}\\[0.5\baselineskip]
\renewcommand\arraystretch{1.4}
\arrayrulewidth 0.5pt \small
\begin{tabular}{cccc}
\hline
&$\|d_{w,r}u_h-u'\|_h$&\quad$\|u-u_h^0\|$&\quad$\max|u_h^i-u(x_i)|$\\
mesh $h$&error \quad\quad rate&\quad error \qquad\quad rate&\quad error \qquad\quad rate\\
\hline
1/4&0.0154\quad\ \ \ \ \ \  \ \ \ \ - &\quad 0.0009 \qquad\ \ \ - &\quad 0.0003 \qquad \ \ \ \ \ \ \ \ \ -\\
1/8&0.0020\quad\ \ \   2.9797 &\qquad5.5590e-5\qquad 3.9842         &\quad 1.7547e-5\qquad\ \ \ 4.0690\\
1/16&2.4534e-4\quad 2.9952 &\qquad3.4831e-6\qquad 3.9964        &\quad 1.1189e-6\qquad\ \ \ 3.9710\\
1/32&3.0693e-5\quad 2.9988 &\qquad 2.1785e-7\qquad 3.9989       &\quad 6.9728e-8\qquad\ \ \ 4.0043\\
1/64&3.8374e-6\quad 2.9997 &\qquad 1.3651e-8\qquad 3.9963       &\quad 4.3549e-9\qquad\ \ \ 4.0010\\
 \hline
\end{tabular}
\end{center}
\vspace{0.2cm}

\begin{center}
{\small{\bf Table 5.3}\quad History of convergence for $k=2$}\\[0.5\baselineskip]
\renewcommand\arraystretch{1.4}
\arrayrulewidth 0.5pt \small
\begin{tabular}{cccc}
\hline
&$\|d_{w,r}u_h-u'\|_h$&\quad$\|u-u_h^0\|$&\quad$\max|u_h^i-u(x_i)|$\\
mesh $h$&error \quad\quad rate&\quad error \qquad\quad rate&\quad error \qquad\quad rate\\
\hline
1/4&0.0008\quad\ \ \ \ \ \ \ \ \ \ - &\quad 2.0906e-5 \qquad\ \ \ - &\quad 1.1846e-6 \qquad \ \ \ \ \ \ \ \ \ \ -\\
1/8&5.1694e-5\quad 3.9944 &\qquad 6.5039e-7\qquad 5.0065         &\quad 1.7776e-8\qquad\ \  \ \ 6.0583\\
1/16&3.2341e-6\quad 3.9986 &\qquad2.0305e-8\qquad 5.0014        &\quad 2.7789e-10\qquad \ \ 5.9993\\
1/32&2.0214e-7\quad 3.9999 &\qquad 6.3690e-10\qquad 4.9947       &\quad 4.2230e-12\qquad\ \ \ 6.0401\\
1/64&1.2594e-8\quad 4.0045 &\qquad 1.9884e-11\qquad 5.0040       &\quad 6.5939e-14\qquad\ \ \ 6.0001\\
 \hline
\end{tabular}
\end{center}
\section*{Acknowledgments}
This work was supported by the National Natural Science Funds of
China, No. 11371081; and the State Key Laboratory of Synthetical
Automation for Process Industries Fundamental Research Funds, No.
2013ZCX02.


\baselineskip 0.5cm
\begin{thebibliography}{26}
\bibitem{Wang}{J. Wang, X. Ye, A weak Galerkin finite element method for
second-order elliptic problems, J. Comput. Appl. Math. 241 (2013)
103-115.}

\bibitem{Wang2}{J. Wang, X. Ye, A weak Galerkin mixed finite element method for
second-order elliptic problems, Math. Comp. 83 (2014) 2101¨C2126.}

\bibitem{Lin}{L. Mu, J. Wang, Y. Wang, X. Ye, A computational study of the
weak Galerkin method for second order elliptic equations, Numer.
Algor. 63 (2012) 753-777.}

\bibitem{Lin1}{L. Mu, J. Wang, G. Wei, X. Ye, S. Zhao, Weak Galerkin methods
for second order elliptic interface problems, J. Comput. Phys. 250
(2013) 106-125.}

\bibitem{Lin2}{L. Mu, J. Wang, X. Ye, Weak Galerkin finite element methods for
the biharmonic equation on polytopal meshes, Numer. Meth. PDEs. 30
(2014) 1003-1029}

\bibitem{Lin3}{L. Mu, J. Wang, X. Ye, A weak Galerkin finite element
method with polynomial reduction, J. Comp. Appl. Math. 285 (2015),
45-58.}

\bibitem{Chen}{L. Chen, J. Wang, X. Ye, A posteriori error estimates for weak
Galerkin finite element methods for second order elliptic
problems, J. Sci. Comput. 6 (2014) 496-511}

\bibitem{Li}{Q. H. Li, J. Wang, Weak Galerkin finite element methods
for parabolic equations, Numer. Meth. PDEs. 29 (2013) 2004-2024.}

\bibitem{Har}{A. Harris, S. Harris, Superconvergence of weak Galerkin
finite element approximation for second order elliptic problems by
$L_2$-projections, Appl. Math. Comp., 227 (2014) 610-621.}

\bibitem{Lin4}{L. Mu, X. Wang, X. Ye, A modified weak Galerkin finite
element method for the Stokes equations, J. Comp. Appl. Math. 275
(2015) 79-90.}

\bibitem{Gao}{F. Gao, X. Wang, A modified weak Galerkin finite
element method for a class of parabolic problems, J. Comp. Appl.
Math. 271 (2014) 1-19.}

\bibitem{Wang0}{X. Wang, N.S. Malluwawadu, F. Gao, T.C. McMillan, A modified
weak Galerkin finite element method, J. Comp. Appl. Math. 271
(2014) 319-327.}

\bibitem{Yang}{M. Yang, Couplings of mixed finite element and weak Galerkin
methods for elliptic problems, J. Appl. Math. Comput. 47 (2015)
327-343.}

\bibitem{Arn}{D. Arnold, F. Brezzi, B. Cockburn, D. Marini, Unified analysis
of discontinuous Galerkin methods for elliptic problems, SIAM J.
Numer. Anal. 39 (2012) 1749-1779.}

\bibitem{Cock}{B. Cockburn, J. Gopalakrishnan, R. Lazarov, Unified hybridization
of discontinuous Galerkin, mixed and continuous Galerkin methods
for second-order elliptic problems, SIAM J. Numer. Anal. 47 (2009)
1319-1365.}

\bibitem{Zhang}{T. Zhang, Theory and Method for Discontinuous
Finite Element, Science Press, Beijing, 2012}


\end{thebibliography}
\end{document}